\documentclass[a4paper,oneside,reqno,12pt]{report}
\usepackage{amsmath,amssymb}
\usepackage{texdraw}
\usepackage{amscd,amsfonts}
\usepackage{amsthm}
\usepackage[ansinew]{inputenc}
\usepackage{amsbsy}
\usepackage[french]{babel}

\usepackage[dvips]{color}
\usepackage{graphics}
\usepackage{graphicx}

\newtheorem{thm}{Theorem}
\numberwithin{thm}{section}

\theoremstyle{definition}

\newtheorem{dfn}{Definition}[section]
\newtheorem{rl}[thm]{Rule}

    \DeclareMathOperator{\cd}{\cd}

\newcommand{\dl}{\vskip.5cm}
\newcommand{\mrg}{\noindent}
\newcommand{\pr}{\dl\mrg\textit{ Proof. }}

\begin{document}

\noindent\textbf{\Large EXPRESSION OF A TENSOR COMMUTATION MATRIX IN TERMS OF THE GENERALIZED GELL-MANN MATRICES}\\

\noindent RAKOTONIRINA Christian\\
\noindent {\emph{Institut Supérieur de Technologie d'Antananarivo,
IST-T, BP 8122,\\ Madagascar}}\\
E-mail:rakotopierre@refer.mg

\begin{center}
\textbf{Abstract}
\end{center} We have expressed the tensor commutation matrix $n\otimes n$ as linear combination of the tensor products of
the generalized Gell-Mann matrices. The tensor commutation matrices
$3\otimes2$ and $2\otimes3$ have been expressed in terms of the
classical Gell-Mann matrices and the Pauli matrices.
\section*{Introduction}

When we had worked on RAOELINA ANDRIAMBOLOLONA idea on the using
tensor product in Dirac equation \cite{rak03}, \cite{wan01} we had
met the unitary
matrix\\
\begin{center}
$U_{2\otimes2}\;=\;\left(
\begin{array}{cccc}
  1 & 0 & 0 & 0 \\
  0 & 0 & 1 & 0 \\
  0 & 1 & 0 & 0 \\
  0 & 0 & 0 & 1 \\
\end{array}
\right)$
\end{center}
This matrix is frequently found in quantum information theory
\cite{fuj01}, \cite{fad95}, \cite{fra02} where one write, by using
the Pauli matrices \cite{fuj01},\cite{fad95},\cite{fra02}
\begin{center}\begin{equation*}U_{2\otimes2}\;=\;\frac{1}{2}\;I_2\otimes I_2
+ \frac{1}{2}\sum_{i=1}^{3}\sigma_i\otimes\sigma_i\;
\end{equation*}\end{center}
with
 $I_2$ the $2\times2$ unit matrix. We call this matrix a tensor
commutation matrix $2\otimes2$. The tensor commutation matrix
$3\otimes3$ is expressed by using the Gell-Mann matrices under the
following form \cite{raksub}
\begin{center}\begin{equation*}
U_{3\otimes3}\;=\;\frac{1}{3}\;I_3\otimes I_3 +
\frac{1}{2}\sum_{i=1}^{8}\lambda_i\otimes\lambda_i\;
\end{equation*}\end{center}

We have to talk a bit about different types of matrices because in
the generalization of the above formulas we will consider the
commutation matrix  as a matrix of fourth order tensor and in
expressing the commutation matrices $U_{3\otimes2}$,
$U_{2\otimes3}$, at the last section, a commutation matrix will be
considered as matrix of second order tensor.\\
 $\mathcal{M}_{m\times
n}\left(\mathbb{C}\right)$ denotes the set of $m\times n$ matrices
whose elements are complex numbers.
\section{Tensor product of matrices}
\subsection{Matrices}
If the elements of a matrix are considered as the components of a
second order tensor, we adopt the habitual notation for a matrix,
without parentheses  inside, whereas if the elements of the matrix
are, for instance, considered as the components of sixth order
tensor, three times covariant and three times contravariant,  then
we represent the matrix of the following way, for example
\begin{center}
 $M\;=\;\left(
\begin{array}{cc}
  \left(
\begin{array}{cc}
  \left(
\begin{array}{cc}
  1 & 0 \\
  1 & 1 \\
\end{array}
\right) & \left(
\begin{array}{cc}
  1 & 1 \\
  3 & 2 \\
\end{array}
\right) \\
  \left(
\begin{array}{cc}
  0 & 0 \\
  0 & 0 \\
\end{array}
\right) & \left(
\begin{array}{cc}
  1 & 1 \\
  1 & 1 \\
\end{array}
\right) \\
\end{array}
\right) & \left(
\begin{array}{cc}
  \left(
\begin{array}{cc}
  1 & 0 \\
  1 & 2 \\
\end{array}
\right) & \left(
\begin{array}{cc}
  7 & 8 \\
  9 & 0 \\
\end{array}
\right) \\
  \left(
\begin{array}{cc}
  3 & 4 \\
  5 & 6 \\
\end{array}
\right) & \left(
\begin{array}{cc}
  9 & 8 \\
  7 & 6 \\
\end{array}
\right) \\
\end{array}
\right) \\
  \left(
\begin{array}{cc}
  \left(
\begin{array}{cc}
  1 & 1 \\
  1 & 1 \\
\end{array}
\right) & \left(
\begin{array}{cc}
  0 & 0 \\
  3 & 2 \\
\end{array}
\right) \\
  \left(
\begin{array}{cc}
  4 & 5 \\
  1 & 6 \\
\end{array}
\right) & \left(
\begin{array}{cc}
  1 & 7 \\
  8 & 9 \\
\end{array}
\right) \\
\end{array}
\right) & \left(
\begin{array}{cc}
  \left(
\begin{array}{cc}
  5 & 4 \\
  3 & 2 \\
\end{array}
\right) & \left(
\begin{array}{cc}
  1 & 0 \\
  1 & 2 \\
\end{array}
\right) \\
  \left(
\begin{array}{cc}
  3 & 4 \\
  5 & 6 \\
\end{array}
\right) & \left(
\begin{array}{cc}
  7 & 8 \\
  9 & 0 \\
\end{array}
\right) \\
\end{array}
\right) \\
  \left(
\begin{array}{cc}
  \left(
\begin{array}{cc}
  1 & 2 \\
  3 & 4 \\
\end{array}
\right) & \left(
\begin{array}{cc}
  9 & 8 \\
  7 & 6 \\
\end{array}
\right) \\
  \left(
\begin{array}{cc}
  5 & 6 \\
  7 & 8 \\
\end{array}
\right) & \left(
\begin{array}{cc}
  5 & 4 \\
  3 & 2 \\
\end{array}
\right) \\
\end{array}
\right) & \left(
\begin{array}{cc}
  \left(
\begin{array}{cc}
  9 & 8 \\
  7 & 6 \\
\end{array}
\right) & \left(
\begin{array}{cc}
  5 & 4 \\
  3 & 2 \\
\end{array}
\right) \\
  \left(
\begin{array}{cc}
  1 & 0 \\
  1 & 2 \\
\end{array}
\right) & \left(
\begin{array}{cc}
  3 & 4 \\
  5 & 6 \\
\end{array}
\right) \\
\end{array}
\right) \\
\end{array}
\right)$
 \end{center}

$M\;=\;\left(
  M_{j_1j_2j_3}^{i_1i_2i_3}\right)$\\

$i_1i_2i_3$\;=\;$111,\; 112,\; 121,\; 122,\; 211,\; 212,\; 221,\;
222,\; 311,\;312,\;
321,\; 322$\\
row indices\\

 $j_1j_2j_3$\;=\;$111,\; 112,\; 121,\; 122,\; 211,\; 212,\; 221,\; 222$\\
 column indices\\
 The first indices $i_1$  and $j_1$  are the indices of the outside parenthesis
 which we call the first order parenthesis; the second indices $i_2$  and
 $j_2$
 are the indices of the next parentheses  which we call the second
 order parentheses; the third indices $i_3$  and $j_3$  are the indices of the
 most interior parentheses, of this example, which we call third
 order parentheses. So, for instance, $M_{121}^{321}$\;=\;5. \\
 If we delete the third order parenthesis, then the elements of the matrix $M$
  are considered as the components of a forth order tensor, twice
  contravariant and twice covariant.\\
A matrix is a diagonal matrix if deleting the interior parentheses
we have a habitual diagonal matrix.\\
 A matrix is a symmetric (resp.antisymmetric) matrix if deleting the interior parentheses we have a
habitual symmetric (resp. antisymmetric) matrix.\\
 We identify one matrix to another matrix if after deleting the interior parentheses they
are the same matrix.
\subsection{Tensor product of matrices}
\begin{dfn}
Consider $A=\left(A^i_j\right)\in\mathcal{M}_{m\times
n}\left(\mathbb{K}\right)$,
$B=\left(B^i_j\right)\in\mathcal{M}_{p\times
r}\left(\mathbb{K}\right)$. The matrix defined by
\begin{center}
$A\otimes B\;=\;\left(
\begin{array}{ccccc}
  A^1_1B & \ldots & A^1_jB & \ldots & A^1_nB \\
  \vdots &  & \vdots &  & \vdots \\
  A^i_1B & \ldots & A^i_jB & \ldots & A^i_nB \\
  \vdots &  & \vdots &  & \vdots \\
  A^m_1B & \ldots & A^m_jB & \ldots & A^m_nB \\
\end{array}%
\right)$
\end{center}
is called the tensor product of the matrix $A$ by the matrix $B$.
\begin{center}
$A\otimes B\in\mathcal{M}_{mp\times nr}\left(\mathbb{K}\right)$
\end{center}
\begin{center}
$A \otimes
B\;=\;\left(C_{j_1j_2}^{i_1i_2}\right)\;=\;\left(A^{i_1}_{j_1}B^{i_2}_{j_2}\right)$
\end{center}
(cf. for example \cite{fuj01} )
 where,\\
  $i_1i_2$ are row indices\\
  $j_1j_2$ are column indices.\\
\end{dfn}
\section{Generalized Gell-Mann matrices}
Let us fix $n\in\mathbb{N}$, $n\geq2$ for all continuation. The
generalized Gell-Mann matrices or $n\times n$-Gell-Mann matrices are
the traceless hermitian $n\times n$ matrices $\Lambda_{1}$,
$\Lambda_{2}$, \ldots, $\Lambda_{n^{2}-1}$ which satisfy the
relation $Tr\left(\Lambda_{i}, \Lambda_{j}\right)\;=\;2\delta_{ij}$,
for all $i$, $j\in \{1, 2, \ldots, n^{2}-1\}$, where
$\delta_{ij}\;=\;
\delta^{ij}\;=\;\delta^{i}_{j}$ the Kronecker symbol \cite{nar89}.\\
However, for the demonstration of the Theorem \ref{thm32} below,
denote, for $1\leq i < j\leq n$, the
$C_{n}^{2}\;=\;\frac{n!}{2!(n-2)!}$ $n\times n$-Gell-Mann matrices
which are symmetric with all elements 0 except the $i$-th row $j$-th
column and the $j$-th row  $i$-th column which are equal to 1, by
$\Lambda^{(ij)}$; the $C_{n}^{2}\;=\;\frac{n!}{2!(n-2)!}$ $n\times
n$-Gell-Mann matrices which are antisymmetric with all elements are
0 except  the $i$-th row $j$-th column which is equal $-i$ and the
$j$-th row $i$-th column which is equal to $i$  , by
$\Lambda^{[ij]}$ and by $\Lambda^{(d)}$,$1 \leq d \leq n-1$, the
following $(n - 1)$
$n\times n$-Gell-Mann matrices which are diagonal:\\

$\Lambda^{(1)}\;=\;\left(
                     \begin{array}{cccccc}
                       1 & 0 &  & \ldots &  & 0 \\
                       0 & -1 &  &  &  &  \\
                        &  & 0 &  &  & \vdots \\
                       \vdots &  &  & \ddots &  &  \\
                        &  &  &  & \ddots &  \\
                       0 &  & \ldots &  &  & 0 \\
                     \end{array}
                   \right)$, $\Lambda^{(2)}\;=\;\frac{1}{\sqrt{3}}\left(
                               \begin{array}{cccccc}
                                 1 & 0 &  & \ldots &  & 0 \\
                                 0 & 1 &  &  &  &  \\
                                  &  & -2 &  &  & \vdots \\
                                 \vdots &  &  & 0 &  &  \\
                                  &  &  &  & \ddots &  \\
                                 0 &  & \ldots &  &  & 0 \\
                               \end{array}
                             \right)$,\\

                             \ldots,

                             $\Lambda^{(n-1)}\;=\;\frac{1}{\sqrt{C_{n}^{2}}}\left(
                                                                              \begin{array}{cccccc}
                                                                                1 & 0 &  & \ldots &  & 0 \\
                                                                                0 & 1 &  &  &  &  \\
                                                                                 &  & 1 &  &  & \vdots \\
                                                                                \vdots &  &  & \ddots &  &  \\
                                                                                 &  &  &  & 1 &  \\
                                                                                0 &  & \ldots &  &  & -(n-1) \\
                                                                              \end{array}
                                                                            \right)$\\

                                      For $n\;=\;2$ we have the Pauli matrices.

\section{Tensor commutation matrices}
For $p$, $q\in\mathbb{N}$, $p\geq 2$, $q\geq 2$, we call tensor
commutation matrices $p\otimes q$ the permutation matrix
$U_{p\otimes q}\in \mathcal{M}_{pq\times pq}\left(\mathbb{C}\right)$
formed by 0 and 1, verifying the property
\begin{center}
$U_{p\otimes q}.(a\otimes b) \;=\; b\otimes a$
\end{center}
 for all $a\in
\mathcal{M}_{p\times 1}\left(\mathbb{C}\right)$, $b\in
\mathcal{M}_{q\times 1}\left(\mathbb{C}\right)$.\\

 Considering $U_{p\otimes q}$ as a matrix of a second order tensor,
 we can construct it by using
the following rule \cite{raksub}.
 \numberwithin{thm}{section}
\begin{rl}\label{thm31}
 Let us start in putting 1 at first
row and first column, after that let us pass into second column in
going down at the rate of $p$ rows and  put 1 at this place,  then
pass into third column in going down at the rate of $p$ rows and put
1,and so on until there is only for us $p-1$ rows for going down
(then we have obtained as number of 1 : $q$. Then pass into the next
column which is the $(q + 1)$-th column, put 1 at the second row of
this column and repeat the process until we have only $p-2$  rows
for going down (then we have obtained as number of 1 : $2q$). After
that pass into the next column which is the $(2q + 2)$ -th column,
put 1 at the third row of this column and repeat the process until
we have only $p - 3$ rows for going down (then we have obtained as
number of 1 : $3q$). Continuing in this way we will have that the
element at $p\times q$-th row and $p\times q$-th column is 1.  The
other elements are 0.
\end{rl}

\begin{thm}\label{thm32}We have

\begin{center}
$U_{n\otimes n}\;=\;\frac{1}{n}I_{n}\otimes
I_{n}+\frac{1}{2}\displaystyle\sum_{i=1}^{n^{2}-1}\Lambda_{i}\otimes\Lambda_{i}$
\end{center}
\end{thm}
\pr
\begin{center}
$I_{n}\otimes
I_{n}\;=\;\left(\delta_{j_{1}j_{2}}^{i_{1}i_{2}}\right)\;=\;
\left(\delta_{j_{1}}^{i_{1}}\delta_{j_{2}}^{i_{2}}\right)$
\end{center}
\numberwithin{equation}{section}
\begin{equation}\label{e31}
U_{n\otimes n}\;=\;
\left(\delta_{j_{2}}^{i_{1}}\delta_{j_{1}}^{i_{2}}\right)
\end{equation}
where,\\
  $i_1i_2$ are row indices\\
  $j_1j_2$ are column indices\cite{fuj01}.\\
  Consider at first, the $C_{n}^{2}$ symmetric $n\times n$-Gell-Mann
  matrices which can be written\\
  \begin{equation*}
  \begin{split}
  \Lambda^{(ij)}\;&=
  \;\left({\Lambda^{(ij)}}_{k}^{l}\right)_{1\leq l\leq n, 1\leq k\leq n}
  \;\\
  &=  \;\left(\delta^{il}\delta^{j}_{k}\right)_{1\leq l\leq n, 1\leq k\leq n}+
  \left(\delta^{jl}\delta^{i}_{k}\right)_{1\leq l\leq n, 1\leq k\leq
  n}\;\\
  &=\;\left(\delta^{il}\delta^{j}_{k}+
  \delta^{jl}\delta^{i}_{k}\right)_{1\leq l\leq n, 1\leq k\leq n}
  \end{split}
  \end{equation*}
  Then
  \begin{equation*}
\Lambda^{(ij)}\otimes
\Lambda^{(ij)}\;=\;\left({\left(\Lambda^{(ij)}\otimes
\Lambda^{(ij)}\right)}_{k_{1}k_{2}}^{l_{1}l_{2}}\right)\;=\;\left(\delta^{il_{1}}\delta^{j}_{k_{1}}+
  \delta^{jl_{1}}\delta^{i}_{k_{1}}\right)\left(\delta^{il_{2}}\delta^{j}_{k_{2}}+
  \delta^{jl_{2}}\delta^{i}_{k_{2}}\right)
\end{equation*}
$l_{1}l_{2}$ row indices, $k_{1}k_{2}$ column indices.\\
That is
\begin{equation*}
 {\left(\Lambda^{(ij)}\otimes
\Lambda^{(ij)}\right)}_{k_{1}k_{2}}^{l_{1}l_{2}}\;=\;\delta^{il_{1}}\delta^{j}_{k_{1}}
\delta^{il_{2}}\delta^{j}_{k_{2}}+\delta^{il_{1}}\delta^{j}_{k_{1}}\delta^{jl_{2}}\delta^{i}_{k_{2}}
+\delta^{jl_{1}}\delta^{i}_{k_{1}}\delta^{il_{2}}\delta^{j}_{k_{2}}
+\delta^{jl_{1}}\delta^{i}_{k_{1}}\delta^{jl_{2}}\delta^{i}_{k_{2}}
\end{equation*}
  The $C_{n}^{2}$ antisymmetric $n\times n$-Gell-Mann
  matrices can be written\\
\begin{equation*}
\Lambda^{[ij]}\;=
  \;\left({\Lambda^{[ij]}}_{k}^{l}\right)_{1\leq l\leq n, 1\leq k\leq n}
  \;=\;\left(-i\delta^{il}\delta^{j}_{k}+
  i\delta^{jl}\delta^{i}_{k}\right)_{1\leq l\leq n, 1\leq k\leq n}
\end{equation*}
 Then
  \begin{equation*}
\Lambda^{[ij]}\otimes
\Lambda^{[ij]}\;=\;\left({\left(\Lambda^{[ij]}\otimes
\Lambda^{[ij]}\right)}_{k_{1}k_{2}}^{l_{1}l_{2}}\right)
\end{equation*}
\begin{equation*}
{\left(\Lambda^{[ij]}\otimes
\Lambda^{[ij]}\right)}_{k_{1}k_{2}}^{l_{1}l_{2}}\;=\;-\delta^{il_{1}}\delta^{j}_{k_{1}}
\delta^{il_{2}}\delta^{j}_{k_{2}}+\delta^{il_{1}}\delta^{j}_{k_{1}}\delta^{jl_{2}}\delta^{i}_{k_{2}}
+\delta^{jl_{1}}\delta^{i}_{k_{1}}\delta^{il_{2}}\delta^{j}_{k_{2}}
-\delta^{jl_{1}}\delta^{i}_{k_{1}}\delta^{jl_{2}}\delta^{i}_{k_{2}}
\end{equation*}
and
\begin{multline*}
  \displaystyle\sum_{1\leq i< j\leq n}{\left(\Lambda^{(ij)}\otimes
\Lambda^{(ij)}\right)}_{k_{1}k_{2}}^{l_{1}l_{2}}+\displaystyle\sum_{1\leq
i< j\leq n}{\left(\Lambda^{[ij]}\otimes
\Lambda^{[ij]}\right)}_{k_{1}k_{2}}^{l_{1}l_{2}}\\
=\;2\displaystyle\sum_{1\leq i< j\leq
n}\left(\delta^{il_{1}}\delta^{j}_{k_{1}}\delta^{jl_{2}}\delta^{i}_{k_{2}}
+\delta^{jl_{1}}\delta^{i}_{k_{1}}\delta^{il_{2}}\delta^{j}_{k_{2}}\right)\;\\
=\;2\displaystyle\sum_{i\neq
j}\delta^{il_{1}}\delta^{j}_{k_{1}}\delta^{jl_{2}}\delta^{i}_{k_{2}}
  \end{multline*}
  the $l_{1}l_{2}$-th row, $k_{1}k_{2}$-th column of the matrix\\
   $\displaystyle\sum_{1\leq i< j\leq
  n}\Lambda^{(ij)}\otimes
\Lambda^{(ij)}+\displaystyle\sum_{1\leq i< j\leq
  n}\Lambda^{[ij]}\otimes
\Lambda^{[ij]}$.\\

Now, consider the diagonal $n\times n$-Gell-Mann matrices. Let
$d\in\mathbb{N}$, $1 \leq d \leq n-1$,
\begin{center}
$\Lambda^{(d)}\;=\;\frac{1}{\sqrt{C_{d+1}^{2}}}\left(\delta_{k}^{l}
\displaystyle\sum_{p=1}^{d}\delta_{k}^{p}-d\delta_{k}^{l}\delta_{k}^{d+1}\right)$
\end{center}
and the $l_{1}l_{2}$-th row, $k_{1}k_{2}$-th of the matrix
$\Lambda^{(d)}\otimes \Lambda^{(d)}$ is\\
\begin{equation*}
\begin{split}
\left(\Lambda^{(d)}\otimes
\Lambda^{(d)}\right)_{k_{1}k_{2}}^{l_{1}l_{2}}
&=\frac{1}{C_{d+1}^{2}}\delta_{k_{1}}^{l_{1}}\delta_{k_{2}}^{l_{2}}
\left(\displaystyle\sum_{q=1}^{d}\displaystyle\sum_{p=1}^{d}\delta_{k_{1}}^{q}\delta_{k_{2}}^{p}\right)\\
&\quad-\frac{1}{C_{d+1}^{2}}\delta_{k_{1}}^{l_{1}}\delta_{k_{2}}^{l_{2}}
\left(d\delta_{k_{2}}^{d+1}\displaystyle\sum_{p=1}^{d}\delta_{k_{1}}^{p}\right)\\
&\quad-\frac{1}{C_{d+1}^{2}}\delta_{k_{1}}^{l_{1}}\delta_{k_{2}}^{l_{2}}
\left(d\delta_{k_{1}}^{d+1}\displaystyle\sum_{p=1}^{d}\delta_{k_{2}}^{p}\right)\\
&\quad+\frac{1}{C_{d+1}^{2}}\delta_{k_{1}}^{l_{1}}\delta_{k_{2}}^{l_{2}}
\left(d^{2}\delta_{k_{1}}^{d+1}\delta_{k_{2}}^{d+1}\right)
\end{split}
\end{equation*}
$\Lambda^{(d)}\otimes \Lambda^{(d)}$ is a diagonal matrix, so all
that we have to do is to calculate the elements on the diagonal
where $l_{1} \;=\; k_{1}$ and $l_{2} \;=\; k_{2}$. Then,
\begin{equation*}
\begin{split}
\sum_{d=1}^{n-1}\left(\Lambda^{(d)}\otimes
\Lambda^{(d)}\right)_{k_{1}k_{2}}^{l_{1}l_{2}}
&=\sum_{d=1}^{n-1}\frac{1}{C_{d+1}^{2}}\left(\displaystyle\sum_{q=1}^{d}\delta_{k_{1}}^{q}\right)
\left(\displaystyle\sum_{p=1}^{d}\delta_{k_{2}}^{p}\right)\\
&\quad-\sum_{d=1}^{n-1}\frac{1}{C_{d+1}^{2}}d\delta_{k_{2}}^{d+1}\displaystyle\sum_{p=1}^{d}\delta_{k_{1}}^{p}\\
&\quad-\sum_{d=1}^{n-1}\frac{1}{C_{d+1}^{2}}d\delta_{k_{1}}^{d+1}\displaystyle\sum_{p=1}^{d}\delta_{k_{2}}^{p}\\
&\quad+\sum_{d=1}^{n-1}\frac{1}{C_{d+1}^{2}}d^{2}\delta_{k_{1}}^{d+1}\delta_{k_{2}}^{d+1}
\end{split}
\end{equation*}
the $l_{1}l_{2}$-th row, $k_{1}k_{2}$-th column of the diagonal
matrix $\displaystyle\sum_{d=1}^{n-1}\Lambda^{(d)}\otimes
\Lambda^{(d)}$ with $l_{1}
\;=\; k_{1}$ and $l_{2} \;=\; k_{2}$.\\
Let us distinguish two cases.\\
 \underline{$1^{st}$ case}: $k_{1}\neq1$ or $k_{2}\neq1$\\
 \indent case1: $k_{1}\neq k_{2}$

\text{If $k_{1}<k_{2}$},
\begin{equation*}
\begin{split}
\sum_{d=1}^{n-1}\left(\Lambda^{(d)}\otimes
\Lambda^{(d)}\right)_{k_{1}k_{2}}^{l_{1}l_{2}}\;&=\;\sum_{d=k_{2}}^{n-1}\frac{1}{C_{d+1}^{2}}-\frac{k_{2}-1}{C_{k_{2}}^{2}}\;\\
&\quad=\;2\left[\sum_{d=k_{2}}^{n-1}\left(\frac{1}{d}-\frac{1}{d+1}\right)-\frac{1}{k_{2}}\right]\;\\
&\quad=\;-\frac{2}{n}
\end{split}
\end{equation*}
 \text{Similarly,if $k_{1}>k_{2}$},
$\displaystyle\sum_{d=1}^{n-1}\left(\Lambda^{(d)}\otimes
\Lambda^{(d)}\right)_{k_{1}k_{2}}^{l_{1}l_{2}}\; =\;-\frac{2}{n}$\\
\indent case2: $k_{1}\;=\; k_{2}\;\neq\;1$
\begin{equation*}
\begin{split}
\sum_{d=1}^{n-1}\left(\Lambda^{(d)}\otimes
\Lambda^{(d)}\right)_{k_{1}k_{2}}^{l_{1}l_{2}}\;
&=\;\sum_{d=k_{2}}^{n-1}\frac{1}{C_{d+1}^{2}}+\frac{\left(k_{2}-1\right)^{2}}{C_{k_{2}}^{2}}\;\\
&\quad=\;\frac{2}{k_{2}}-\frac{2}{n}+\frac{\left(k_{2}-1\right)^{2}}{C_{k_{2}}^{2}}\;\\
 &\quad=\;2-\frac{2}{n}
\end{split}
\end{equation*}
\underline{$2^{nd}$ case}: $k_{1}= k_{2}=1$
\begin{equation*}
\sum_{d=1}^{n-1}\left(\Lambda^{(d)}\otimes
\Lambda^{(d)}\right)_{k_{1}k_{2}}^{l_{1}l_{2}}\;=\;\sum_{d=1}^{n-1}\frac{1}{C_{d+1}^{2}}\;=\;2-\frac{2}{n}
\end{equation*}
We can condense these cases in one formula
\begin{equation*}
\sum_{d=1}^{n-1}\left(\Lambda^{(d)}\otimes
\Lambda^{(d)}\right)_{k_{1}k_{2}}^{l_{1}l_{2}}\;
=\;-\frac{2}{n}\delta_{k_{1}}^{l_{1}}\delta_{k_{2}}^{l_{2}}
+2\sum_{i=1}^{n}\delta^{il_{1}}\delta_{k_{1}}^{i}\delta^{il_{2}}\delta_{k_{2}}^{i}
\end{equation*}
which yields the diagonal of the diagonal matrix
$\displaystyle\sum_{d=1}^{n-1}\Lambda^{(d)}\otimes \Lambda^{(d)}$.\\

For all the $n\times n$- Gell-Mann matrices we have\\
\begin{multline*}
\sum_{1\leq i< j\leq n}{\left(\Lambda^{(ij)}\otimes
\Lambda^{(ij)}\right)}_{k_{1}k_{2}}^{l_{1}l_{2}}+\sum_{1\leq i<
j\leq n}{\left(\Lambda^{[ij]}\otimes
\Lambda^{[ij]}\right)}_{k_{1}k_{2}}^{l_{1}l_{2}}+\sum_{d=1}^{n-1}\left(\Lambda^{(d)}\otimes
\Lambda^{(d)}\right)_{k_{1}k_{2}}^{l_{1}l_{2}}\;\\
=\;-\frac{2}{n}\delta_{k_{1}}^{l_{1}}\delta_{k_{2}}^{l_{2}}
+2\sum_{i=1}^{n}\delta^{il_{1}}\delta_{k_{1}}^{i}\delta^{il_{2}}\delta_{k_{2}}^{i}
+2\sum_{i\neq
j}\delta^{il_{1}}\delta_{k_{1}}^{j}\delta^{jl_{2}}\delta_{k_{2}}^{i}\\
=\;-\frac{2}{n}\delta_{k_{1}}^{l_{1}}\delta_{k_{2}}^{l_{2}}
+2\sum_{j=1}^{n}\sum_{i=1}^{n}\delta^{il_{1}}\delta_{k_{1}}^{j}\delta^{jl_{2}}\delta_{k_{2}}^{i}\\
=\;-\frac{2}{n}\delta_{k_{1}}^{l_{1}}\delta_{k_{2}}^{l_{2}}+2\delta_{k_{2}}^{l_{1}}\delta_{k_{1}}^{l_{2}}
\end{multline*}
for all $l_{1}$, $l_{2}$, $k_{1}$,$k_{2}\in\{1, 2, \ldots, n\}$.\\
Hence, by using \eqref{e31}
\begin{equation*}
\displaystyle\sum_{i=1}^{n^{2}-1}\Lambda_{i}\otimes\Lambda_{i}\;=\;-\frac{2}{n}I_{n}\otimes
I_{n}+2U_{n\otimes n}
\end{equation*}
and the theorem is proved.
 \qed
\section{Expression of $U_{3\otimes 2}$ and $U_{2\otimes 3}$ }
In this section we derive formulas for $U_{3\otimes 2}$ and
$U_{2\otimes 3}$, naturally in terms of the Pauli matrices\\
\begin{center}$\sigma_1$ = $\left(
\begin{array}{cc}
  0 & 1 \\
  1 & 0 \\
\end{array}
\right)$, $\sigma_2$\;=\;$\left(
\begin{array}{cc}
  0 & -i \\
  i & 0 \\
\end{array}
\right)$, $\sigma_3$\;=\;$\left(
\begin{array}{cc}
  1 & 0 \\
  0 & -1 \\
\end{array}
\right)$
\end{center}
and the Gell-Mann matrices\\

$\lambda_{1}=\left(
                \begin{array}{ccc}
                  0 & 1 & 0 \\
                  1 & 0 & 0 \\
                  0 & 0 & 0 \\
                \end{array}
              \right)$,\; $\lambda_{2}=\left(
                \begin{array}{ccc}
                  0 & -i & 0 \\
                  i & 0 & 0 \\
                  0 & 0 & 0 \\
                \end{array}
              \right)$,\;
               $\lambda_{3}=\left(
                \begin{array}{ccc}
                  1 & 0 & 0 \\
                  0 & -1 & 0 \\
                  0 & 0 & 0 \\
                \end{array}
              \right)$,\;

               $\lambda_{4}=\left(
                \begin{array}{ccc}
                  0 & 0 & 1 \\
                  0 & 0 & 0 \\
                  1 & 0 & 0 \\
                \end{array}
              \right)$,
               $\lambda_{5}=\left(
                \begin{array}{ccc}
                  0 & 0 & -i \\
                  0 & 0 & 0 \\
                  i & 0 & 0 \\
                \end{array}
              \right)$,\;\\

               $\lambda_{6}=\left(
                \begin{array}{ccc}
                  0 & 0 & 0 \\
                  0 & 0 & 1 \\
                  0 & 1 & 0 \\
                \end{array}
              \right)$,\; $\lambda_{7}=\left(
                \begin{array}{ccc}
                  0 & 0 & 0 \\
                  0 & 0 & -i \\
                  0 & i & 0 \\
                \end{array}
              \right)$,\;\\

               $\lambda_{8}=\frac{1}{\sqrt{3}}\left(
                \begin{array}{ccc}
                  1 & 0 & 0 \\
                  0 & 1 & 0 \\
                  0 & 0 & -2 \\
                \end{array}
              \right)$\\

For $r\in \mathbb{N}^{*}$,  define $E_{ij}^{(r)}$ the elementary
$r\times r$ matrix whose elements are zeros except the $i$-th row
and $j$-th column which is equal to 1. We construct $U_{3\otimes 2}$
by using the Rule \ref{thm31}, and then we have
\begin{equation*}
U_{3\otimes 2} = E_{11}^{(6)} + E_{23}^{(6)} + E_{35}^{(6)} +
E_{42}^{(6)} + E_{54}^{(6)} + E_{66}^{(6)}
\end{equation*}
 Take
\begin{equation*}
  E_{11}^{(6)} =
E_{11}^{(3)}\otimes E_{11}^{(2)}
\end{equation*}
Let
\begin{equation*}
E_{11}^{(3)} = \alpha_{0}I_{3} + \alpha_{3}\lambda_{3} +
\alpha_{8}\lambda_{8}
\end{equation*}
 with  $\alpha_{0}$,   $\alpha_{3}$,   $\alpha_{8}\in \mathbb{C}$, then
 \begin{center}
 $\alpha_{0} = \frac{1}{3}$ ,   $\alpha_{3} = \frac{1}{2}$ ,   $\alpha_{8}
 =\frac{\sqrt{3}}{6}$
 \end{center}
and
\begin{equation*}
E_{11}^{(3)} = \frac{1}{3}I_{3} + \frac{1}{2}\lambda_{3} +
\frac{\sqrt{3}}{6}\lambda_{8}
\end{equation*}
Let
\begin{equation*}
 E_{11}^{(2)} =  \beta_{0}I_{2}  +   \beta_{3}\sigma_{3}
\end{equation*}
   with $\beta_{0}$, $\beta_{3}\in \mathbb{C}$, then
\begin{center}
 $\beta_{0} =\frac{1}{2} $ ,   $\beta_{3} =\frac{1}{2}$
 \end{center}
and
\begin{equation*}
E_{11}^{(2)} =  \frac{1}{2}I_{2}  +   \frac{1}{2}\sigma_{3}
\end{equation*}
 So we have
\begin{equation*}
 E_{11}^{(6)} = = \left(\frac{1}{3}I_{3} + \frac{1}{2}\lambda_{3} +
\frac{\sqrt{3}}{6}\lambda_{8} \right)\otimes \left( \frac{1}{2}I_{2}
+ \frac{1}{2}\sigma_{3} \right)
\end{equation*}
By the similar way, we have
\begin{align*}
 E_{23}^{(6)} &= \left(\frac{1}{2}\lambda_{1}+\frac{i}{2}\lambda_{2}\right)\otimes
 \left(\frac{1}{2}\sigma_{1}- \frac{i}{2}\sigma_{2}\right)\\
E_{35}^{(6)} &= \left(\frac{1}{2}\lambda_{6}+\frac{i}{2}\lambda_{7}\right)\otimes
\left(\frac{1}{2}I_{2}+\frac{1}{2}\sigma_{3}\right)\\
 E_{42}^{(6)} &= \left( \frac{1}{2}\lambda_{1}-\frac{i}{2}\lambda_{2} \right)\otimes
 \left(\frac{1}{2}I_{2}-\frac{1}{2}\sigma_{3} \right)\\
 E_{54}^{(6)} &= \left(\frac{1}{2}\lambda_{6}-\frac{i}{2}\lambda_{7} \right)\otimes
 \left( \frac{1}{2}\sigma_{1}+ \frac{i}{2}\sigma_{2} \right)\\
  E_{66}^{(6)} &=\left( \frac{1}{3}I_{3}-\frac{\sqrt{3}}{3}\lambda_{8} \right)\otimes
  \left( \frac{1}{2}I_{2}- \frac{1}{2}\sigma_{3} \right)
\end{align*}
Hence
\begin{equation*}
\begin{split}
U_{3\otimes 2} &= \left(\frac{1}{3}I_{3} + \frac{1}{2}\lambda_{3} +
\frac{\sqrt{3}}{6}\lambda_{8} \right)\otimes \left( \frac{1}{2}I_{2}
+ \frac{1}{2}\sigma_{3} \right)\\
& +\left(\frac{1}{2}\lambda_{1}+\frac{i}{2}\lambda_{2}\right)\otimes
 \left(\frac{1}{2}\sigma_{1}- \frac{i}{2}\sigma_{2}\right)\\
&+ \left(\frac{1}{2}\lambda_{6}+\frac{i}{2}\lambda_{7}\right)\otimes
\left(\frac{1}{2}I_{2}+\frac{1}{2}\sigma_{3}\right)\\
 &+ \left( \frac{1}{2}\lambda_{1}-\frac{i}{2}\lambda_{2} \right)\otimes
 \left(\frac{1}{2}I_{2}-\frac{1}{2}\sigma_{3} \right)\\
 &+ \left(\frac{1}{2}\lambda_{6}-\frac{i}{2}\lambda_{7} \right)\otimes
 \left( \frac{1}{2}\sigma_{1}+ \frac{i}{2}\sigma_{2} \right)\\
  &+\left( \frac{1}{3}I_{3}-\frac{\sqrt{3}}{3}\lambda_{8} \right)\otimes
  \left( \frac{1}{2}I_{2}- \frac{1}{2}\sigma_{3} \right)
  \end{split}
\end{equation*}
   From analogous way,
\begin{equation*}
\begin{split}
U_{2\otimes 3} &= \left( \frac{1}{2}I_{2} + \frac{1}{2}\sigma_{3}
\right)\otimes \left(\frac{1}{3}I_{3} + \frac{1}{2}\lambda_{3} +
\frac{\sqrt{3}}{6}\lambda_{8} \right)\\
&+\left(\frac{1}{2}\sigma_{1}+
\frac{i}{2}\sigma_{2}\right)\otimes\left(\frac{1}{2}\lambda_{1}-\frac{i}{2}\lambda_{2}\right)\\
&+\left(\frac{1}{2}I_{2}+\frac{1}{2}\sigma_{3}\right)\otimes
\left(\frac{1}{2}\lambda_{6}-\frac{i}{2}\lambda_{7}\right)\\
&+\left(\frac{1}{2}I_{2}-\frac{1}{2}\sigma_{3} \right)\otimes\left(
\frac{1}{2}\lambda_{1}+\frac{i}{2}\lambda_{2} \right) \\
&+ \left( \frac{1}{2}\sigma_{1}- \frac{i}{2}\sigma_{2}
\right)\otimes\left(\frac{1}{2}\lambda_{6}+\frac{i}{2}\lambda_{7}
\right)\\
  &+ \left( \frac{1}{2}I_{2}- \frac{1}{2}\sigma_{3} \right)\otimes
  \left( \frac{1}{3}I_{3}-\frac{\sqrt{3}}{3}\lambda_{8} \right)
\end{split}
\end{equation*}
 We can develop these formulas in employing the distributivity of the tensor product.

\subsection*{Acknowledgements} The author thanks  the referee of an earlier manuscript
for suggesting the topic. The author would like to thank Victor
Razafinjato, Director of Genie civil Departement of IST-T and
Ratsimbarison Mahasedra for encouragement and for critical reading
the manuscript.

\renewcommand{\bibname}{References}

\end{document}